\documentclass[british]
{elsarticle}
\usepackage{amsmath, amssymb, amsthm}
\usepackage{babel,eucal,url,amssymb,
enumerate,amscd,
}

\usepackage[pagebackref]{hyperref}

\begin{document}


\title
{A connection with parallel totally skew-symmetric torsion 
on a class of almost hypercomplex manifolds 
with Hermitian and anti-Hermitian metrics}

\author{Mancho Manev}
\ead{mmanev@uni-plovdiv.bg}
\ead[url]{http://fmi.uni-plovdiv.bg/manev/}
\author{Kostadin Gribachev}
\ead{costas@uni-plovdiv.bg}
\address{Paisii Hilendarski University of Plovdiv \\
Faculty of Mathematics and Informatics\\
236 Bulgaria Blvd., 4003 Plovdiv, Bulgaria}




\newcommand{\ie}{i.e. }
\newcommand{\X}{\mathfrak{X}}
\newcommand{\W}{\mathcal{W}}
\newcommand{\G}{\mathcal{G}}
\newcommand{\K}{\mathcal{K}}
\newcommand{\N}{\mathbb{N}}
\newcommand{\R}{\mathbb{R}}
\newcommand{\s}{\mathfrak{S}}
\newcommand{\n}{\nabla}
\newcommand{\D}{{\rm d}}
\newcommand{\al}{\alpha}
\newcommand{\bt}{\beta}
\newcommand{\gm}{\gamma}
\newcommand{\ea}{\varepsilon_\alpha}
\newcommand{\eb}{\varepsilon_\beta}
\newcommand{\eg}{\varepsilon_\gamma}
\newcommand{\dd}{\textrm{d}}
\newcommand{\sa}{\sum_{\al=1}^3}
\newcommand{\sbt}{\sum_{\bt=1}^3}
\newcommand{\ee}{\end{equation}}
\newcommand{\be}[1]{\begin{equation}\label{#1}}
\def\bea{\begin{eqnarray}} \def\eea{\end{eqnarray}}
\newcommand{\norm}[1]{\left\Vert#1\right\Vert ^2}
\newcommand{\nJ}[1]{\norm{\nabla J_{#1}}}
\newcommand{\sx}{\mathop{\mathfrak{S}}\limits_{x,y,z}}

\newcommand{\thmref}[1]{The\-o\-rem~\ref{#1}}
\newcommand{\propref}[1]{Pro\-po\-si\-ti\-on~\ref{#1}}
\newcommand{\secref}[1]{\S\ref{#1}}
\newcommand{\lemref}[1]{Lem\-ma~\ref{#1}}
\newcommand{\dfnref}[1]{De\-fi\-ni\-ti\-on~\ref{#1}}
\newcommand{\coref}[1]{Corollary~\ref{#1}}

\renewcommand{\thefootnote}{\fnsymbol{footnote}}


\numberwithin{equation}{section}
\newtheorem{thm}{Theorem}[section]
\newtheorem{lem}[thm]{Lemma}
\newtheorem{prop}[thm]{Proposition}
\newtheorem{cor}[thm]{Corollary}

\theoremstyle{definition}
\newtheorem{defn}{Definition}[section]
\newtheorem{rem}{Remark}[section]

\hyphenation{Her-mi-ti-an ma-ni-fold ah-ler-ian}

\begin{keyword}
almost hypercomplex manifold \sep pseudo-Riemannian metric
 \sep anti-Hermitian metric \sep Norden metric \sep indefinite metric
 \sep neutral metric \sep
natural connection \sep parallel structure \sep totally
skew-symmetric torsion \sep strong HKT-connection \sep weak HKT-connection. %
\MSC[2010]{53C26, 53C15, 53C50, 53C55.}
\end{keyword}


\begin{abstract}
The subject of investigations are the almost hypercomplex
manifolds with Hermitian and anti-Hermitian (Norden) metrics. A
linear connection $D$ is introduced such that the structure of
these manifolds is parallel with respect to $D$ and its torsion is
totally skew-symmetric. The class of the nearly K\"ahler manifolds
with respect to the first almost complex structure is of special
interest. It is proved that $D$ has a $D$-parallel torsion and is
weak if it is not flat. Some curvature properties of these
manifolds are studied.
\end{abstract}

\maketitle



\section*{Introduction}
The linear connections whose torsion is a 3-form, \ie connections
with totally skew-symmetric torsion, are known as \emph{K\"ahler
with torsion} (shortly, \emph{KT}-) \emph{connections} or
\emph{Bismut connections} on an almost Hermitian manifold
\cite{Stro,Bis}. There exists, on any Hermitian manifold,  a
unique KT-connection whose torsion tensor is a 3-form \cite{Gau,
Fri-Iv2}. In \cite{Fri-Iv2} all almost contact metric, almost
Hermitian and $G_2$-structures admitting a connection with totally
skew-symmetric torsion tensor are described. Such a connection on
almost complex manifolds with Norden metric is introduced and
investigated in \cite{Mek2, Mek6, Mek7}. A connection of this type
on almost contact manifolds with B-metric (which is the
corresponding metric of the Norden metric in the odd dimensional
case) is introduced in \cite{Man31}.

A good quaternionic analog of K\"ahler geometry is given by the
\emph{hyper-K\"ahler with torsion} (shortly, \emph{HKT-})
geometry.
An HKT-manifold is a hyper-Hermitian manifold 
admitting a hyper-Hermitian connection with totally
skew-symmetric torsion, \ie for which the three KT-connections
associated to the three Hermitian structures 
coincide. This connection is said to be an \emph{HKT-connection}.
This geometry is introduced in \cite{HoPa} and later studied for
instance in \cite{GranPoon, FiGran, BarDotVerb, BarFi, Sw}. A case
of  particular interest is when the torsion 3-form of such
HKT-connection is closed. In this case the HKT-manifold is called
\emph{strong} otherwise it is called \emph{weak}. The HKT-geometry
is a natural generalization of the hyper-K\"ahler geometry, since
when the torsion is zero the HKT-connection coincides with the
Levi-Civita connection. The HKT-connections have applications in
some branches of theoretical and mathematical physics. For
instance,  these connections appear on supersymmetric sigma models
with Wess-Zumino term \cite{GHR, HoPa, HoPa2} and in supergravity
theory \cite{GiPaSte, PaTe}. Some of the applications of KT- and
HKT-geometries in physics are as follows: Strong KT- and
HKT-geometries have applications in type II string theory and in
two-dimensional supersymmetric sigma models \cite{Stro, GHR,
HoPa3}. Weak KT- and HKT-geometries have applications in
supersymmetric quantum mechanics \cite{GiPaSte, CoPa}. Strong and
weak HKT-geometries have applications in the moduli spaces of
gravitational solitons and black holes \cite{GiPaSte, MiStro,
GuPa}.

In this work\footnote{This work was partially supported by
projects RS09-FMI-003 and IS-M-4/2008 of the Scientific Research
Fund at the University of Plovdiv.} we continue the investigations
on
 an almost hypercomplex manifold $(M,H)$ 
with a metric structure $G$, generated by a pseudo-Riemannian
metric
$g$ of neutral signature \cite{GrMa,GrMaDi}.

As it is known, if $g$ is a Hermitian metric on $(M,H)$, the
derived metric structure $G$ is the known hyper-Hermitian
structure. It consists of the given Hermitian metric $g$ with
respect to the three almost complex structures of $H$ and the
three K\"ahler forms associated with $g$ by $H$ (see, e.g.
\cite{AlMa}).
%

Here, the considered metric structure $G$ has a different type of
compatibility with $H$. The structure $G$ is generated by a
neutral metric $g$ such that the first (resp., the other two) of
the almost complex structures of $H$ acts as an iso\-metry (resp.,
act as anti-isometries) with respect to $g$ in each tangent
fibre. 
%
Suppose that the almost complex structures of $H$ act as
iso\-metries or anti-isometries with respect to the metric, then
the existence of an anti-isometry generates exactly the existence
of one more
anti-isometry and an isometry. %
Thus, $G$ contains the metric $g$ and three (0,2)-tensors
associated by $H$ -- a K\"ahler form and two metrics of the same
type.
The existence of such bilinear forms on an almost hypercomplex
manifold is proved in \cite{GrMa}.
The neutral metric $g$ is Her\-mit\-ian with respect to the first
almost complex structure of $H$ and $g$ is an anti-Hermitian (\ie
 Norden) metric regarding the other two almost complex structures
of $H$. For this reason we call the derived manifold $(M,H,G)$ an
\emph{almost hypercomplex manifold with Hermitian and
anti-Hermitian metrics} or shortly an \emph{almost
$(H,G)$-manifold}.

The geometry of an arbitrary almost $(H,G)$-manifold is the
geometry of the hypercomplex structure $H=\{J_1,J_2,J_3\}$ and the
neutral metric $g$ or equivalently -- the geometry of the metric
structure $G=\{g,g_1,g_2,g_3\}$,
$g_\al(\cdot,\cdot):=g(J_\al\cdot,\cdot)$, $\al=1,2,3$. In this
geometry, there are important so-called \emph{natural connections}
for the $(H,G)$-structure (briefly, the $(H,G)$-connections), \ie
$H$ and $g$ are parallel with respect to such a connection.

If the three almost complex structures of $H$ are parallel with
respect to the Levi-Civita connection $\nabla$ of $g$, then we
call such $(H,G)$-manifolds of K\"ahler type
\emph{pseudo-hyper-K\"ahler manifolds} and we denote their class
by $\K$. Therefore, outside of the class $\K$, the Levi-Civita
connection $\n$ is no longer an $(H,G)$-connection. There exist
countless natural connections on an almost $(H,G)$-manifold in the
general case.

In this paper we construct a natural connection $D$  with totally
skew-symmet\-ric torsion tensor on almost hypercomplex manifolds
with Hermitian and Norden metrics. The presence of almost complex
structures with Hermitian and Norden metrics gives us a reason to
restrict our consideration to those classes where the
corresponding natural connections with totally skew-symmetric
torsion exist. Namely, these are the class $\G_1$ of almost
Hermitian manifolds and the class of $\W_3$ of quasi-K\"ahler
manifolds with Norden metric. A central place in our
investigations takes the $\G_1$'s subclass $\W_1$ of nearly
K\"ahler manifolds -- an important part in the theory of geometric
structures of non-integrable type. The study of these manifolds
was initiated by A.~Gray in 1970's as the concept of weak holonomy
and they have been investigated since then by many authors.

The present paper is organized as follows.
In Section~1 we present some necessary facts about the considered
manifolds.
Section~2 is devoted to the study of the properties of the
considered manifolds of a special class. Namely, the manifolds
which belong to $\G_1$ with respect to $J_1$ in the Gray-Hervella
classification and to $\W_3$ with respect to $J_2$ and $J_3$ in
the Ganchev-Borisov classification. After that, we characterize
the more specialized subclass (denoted by $\W_{133}$) of the above
class, where the manifolds are nearly K\"ahlerian with respect to
$J_1$.
In Section~3 we introduce the notion of a pseudo-HKT-connection
(shortly, pHKT-connection) on an almost $(H,G)$-manifold and
determine the class of the considered manifold, where such a
connection exists. It is this class considered in the first part
of Section~2.
Then, we construct the pHKT-connection $D$ in a special case when
the almost $(H,G)$-manifold belongs to $\W_{133}$.
At the end of Section~3 we show that this connection has a
$D$-parallel torsion and prove that any non-$D$-flat $(M,H,G)$ is
a weak pHKT-manifold.

The point in question in this work is the existence and the
geometric characteristics of the considered manifolds with totally
skew-symmetric torsion. The main result of this paper is that the
known KT-connection on a nearly K\"ahler manifold plays the role
of the pseudo-HKT-connection on the corresponding almost
$(H,G)$-manifold.

\section{Almost hypercomplex manifolds with Hermitian and anti-Her\-mit\-ian metrics}

\subsection{The almost $(H,G)$-manifolds}

Let $(M,H)$ be an almost hypercomplex manifold, \ie $M$ is a
$4n$-dimension\-al differentiable manifold and $H=(J_1,J_2,J_3)$
is a triple of almost complex structures with the properties:
\begin{equation}\label{J123} %
J_\al=J_\bt\circ J_\gm=-J_\gm\circ J_\bt, \qquad
J_\al^2=-I%
\end{equation} %
for all cyclic permutations $(\al, \bt, \gm)$ of $(1,2,3)$ and $I$
denotes the identity (\cite{Boy}, \cite{AlMa}).

The standard structure of $H$ on a $4n$-dimensional vector space
with a basis
$\{X_{4k+1},X_{4k+2},X_{4k+3},X_{4k+4}\}_{k=0,1,\dots,n-1}$ has
the form \cite{So}:
\begin{equation}\label{Jdim4n}
\begin{array}{lll}
J_1X_{4k+1}= X_{4k+2},\;\; & J_2X_{4k+1}=X_{4k+3},\;\; & J_3X_{4k+1}=-X_{4k+4},\\[4pt]
J_1X_{4k+2}=-X_{4k+1},\;\; & J_2X_{4k+2}=X_{4k+4}, \;\; & J_3X_{4k+2}=X_{4k+3}, \\[4pt]
J_1X_{4k+3}=-X_{4k+4},\;\; & J_2X_{4k+3}=-X_{4k+1}, \;\; & J_3X_{4k+3}=-X_{4k+2},\\[4pt]
J_1X_{4k+4}=X_{4k+3},\;\; & J_2X_{4k+4}=-X_{4k+2}, \;\; &
J_3X_{4k+4}=X_{4k+1}.
\end{array}
\end{equation}


Further, 
$x$, $y$, $z$, $w$ 
will stand for arbitrary differentiable vector fields on $M$.

Let $g$ be a pseudo-Riemannian metric on $(M,H)$ with the
properties
\begin{equation}\label{gJJ} %
g(x,y)=\ea g(J_\al x,J_\al y), \end{equation} %
where %
\begin{equation}\label{ea}%
 \ea=
\begin{cases}
\begin{array}{ll}
1, \quad & \al=1;\\
-1, \quad & \al=2;3.
\end{array}
\end{cases}
\end{equation}
In other words, for $\al=1$, the metric $g$ is Hermitian with
respect to $J_1$, and in the case  $\al=2$ or $\al=3$, the metric
$g$ is an anti-Hermitian (\ie Norden) metric with respect to
$J_\al$ ($\al=2$ or $\al=3$, respectively) \cite{GaBo}. Moreover,
the associated bilinear forms $g_1$, $g_2$, $g_3$ are determined
by
\begin{equation}\label{gJ} %
g_\al(x,y)=g(J_\al x,y)=-\ea g(x,J_\al y),\qquad
\al=1,2,3. %
\end{equation}
Following \eqref{gJJ} and \eqref{gJ}, the metric $g$ and the
associated bilinear forms $g_2$ and $g_3$ are necessarily
pseudo-Riemannian metrics of neutral signature $(2n,2n)$. The
associated bilinear form $g_1$ is the associated (K\"ahler)
2-form.

A structure $(H,G)=(J_1,J_2,J_3,g,g_1,g_2,g_3)$ is introduced and
investigated in \cite{GrMa}, \cite{GrMaDi}, \cite{ManSek},
\cite{Man19} and \cite{Man28}. The cases when the original metric
$g$ is Hermitian or anti-Hermitian with respect to the almost
complex structures of $H$ are considered. In \cite{GrMa} it is
proved that the unique possibility for an anti-Hermitian metric to
be considered on an almost hypercomplex manifold is the case when
the given metric is Hermitian with respect to the first and
moreover it is an anti-Hermitian metric with respect to other two
structures of $H$. Therefore, we call $(H,G)$ an \emph{almost
hypercomplex structure with Hermitian and anti-Hermitian
met\-rics} on $M$ (or, in short, an \emph{almost
$(H,G)$-structure} on $M$). Then, we call a manifold with such a
structure briefly an \emph{almost $(H,G)$-manifold}.

The structural tensors of an almost $(H,G)$-manifold are the three
$(0,3)$-tensors determined by
\begin{equation}\label{F}
F_\al (x,y,z)=g\bigl( \left( \n_x J_\al
\right)y,z\bigr)=\bigl(\n_x g_\al\bigr) \left( y,z \right),\qquad
\al=1,2,3,
\end{equation}
where $\n$ is the Levi-Civita connection generated by $g$.
%

The tensors $F_\al$ have the following fundamental identities:
\begin{equation}\label{34}%
\begin{array}{c}
    F_\al(x,y,z)=-\ea F_\al(x,z,y)=-\ea F_\al(x,J_\al y,J_\al
    z),\\[4pt]
    F_\al(x,J_\al y,z)=\ea F_\al(x,y,J_\al z);
\end{array}
\end{equation}
%
\begin{equation}\label{F-prop-1}
\begin{array}{l}
    F_\al(x,y,z)=F_\bt(x,J_\gm y,z)-\eb F_\gm(x,y,J_\bt z)\\[4pt]
    \phantom{F_\al(x,y,z)}=-F_\gm(x,J_\bt y,z)+\eg F_\bt(x,y,J_\gm
    z)
\end{array}
\end{equation}
for all cyclic permutations $(\al, \bt, \gm)$ of $(1,2,3)$.

In \cite{GrMaDi} we study a special class $\K$ of the
$(H,G)$-mani\-folds -- the so-called there
\emph{pseudo-hyper-K\"ahler manifolds}. The manifolds from the
class $\K$ are the $(H,G)$-mani\-folds for which the complex
structures $J_\al$ are parallel with respect to the Levi-Civita
connection $\n$, generated by $g$, for all $\al=1,2,3$.
%

As $g$ is an indefinite metric, there exist isotropic vectors
on $M$, \ie \(g(x,x)=0\) for a nonzero vector \(x\). 
In \cite{GrMa} we define the invariant square norm
\begin{equation}\label{nJ}
\nJ{\alpha}= g^{ij}g^{kl}g\bigl( \left( \nabla_i J_\alpha \right)
e_k, \left( \nabla_j J_\alpha \right) e_l \bigr),
\end{equation}
where $\{e_i\}_{i=1}^{4n}$ is an arbitrary basis of the tangent
space $T_pM$ at an arbitrary point $p\in M$. We say that an almost
$(H,G)$-manifold is an \emph{isotropic pseudo-hyper-K\"ahler
manifold} if $\nJ{\alpha}=0$ for all $J_\alpha$ of $H$. Clearly,
if $(M,H,G)$ is a pseudo-hyper-K\"ahler manifold, then it is an
isotropic pseudo-hyper-K\"ahler manifold. The inverse statement
does not hold. For instance, in \cite{GrMa} we have constructed an
almost $(H,G)$-manifold on a Lie group, which is an isotropic
pseudo-hyper-K\"ahler manifold but it is not a
pseudo-hyper-K\"ahler manifold.

\subsection{Properties of the K\"ahler-like tensors}

A tensor $L$ of type (0,4) with the pro\-per\-ties:%
\begin{equation}\label{curv}%
\begin{split}%
&L(x,y,z,w)=-L(y,x,z,w)=-L(x,y,w,z),\\[4pt]
&L(x,y,z,w)+L(y,z,x,w)+L(z,x,y,w)=0
\end{split}%
\end{equation} %
is called a \emph{curvature-like tensor}. The last equality of
\eqref{curv} is known as the first Bianchi identity of a
curvature-like tensor $L$.

We say that a curvature-like tensor $L$ is a \emph{K\"ahler-like
tensor} on an almost $(H,G)$-manifold when $L$ satisfies the properties: %
\begin{equation}\label{L-kel}%
\begin{array}{l}
L(x,y,z,w)=\ea L(x,y,J_\al z,J_\al w)=\ea L(J_\al x,J_\al y,z,w),
\end{array}
\end{equation} where $\ea$ is determined by \eqref{ea}.

Let the curvature tensor $R$ of the Levi-Civita connection
$\nabla$, generated by $g$, be defined, as usual, by
$R(x,y)z=\nabla_x \nabla_y z - \nabla_y \nabla_x z -
\nabla_{\left[x,y\right]} z$. The corresponding $(0,4)$-tensor,
denoted by the same letter, is determined by
$R(x,y,z,w)=g\left(R(x,y)z,w\right)$. Obviously, $R$ is a
K\"ahler-like tensor on an arbitrary pseudo-hyper-K\"ahler
manifold.

A K\"ahler-like tensor $L$ on an arbitrary almost $(H,G)$-manifold
has the same properties \eqref{curv} and \eqref{L-kel} of $R$ on a
pseudo-hyper-K\"ahler manifold. Thus, we obtain the following
geometric characteristic of the K\"ahler-like tensors on an almost
$(H,G)$-manifold, similarly to Theorem~2.3 in \cite{GrMaDi}, where
it is proved that the hyper-K\"ahler $(H,G)$-manifolds are flat,
\ie $R=0$ in $\K$. Using the same idea, we establish the
truthfulness of the following
\begin{prop}[\cite{Man28}]\label{th-0}
Every K\"ahler-like tensor on an almost $(H,G)$-mani\-fold is
ze\-ro.  \phantom{} \hfill $\Box$
\end{prop}

\section{Properties of the $(H,G)$-manifolds of a certain class}

\subsection{Class $\G_1(J_1)\cap\W_3(J_2)\cap\W_3(J_3)$}

For the sake of our further purposes, we restrict the class of
$(M,H,G)$ to the class $\G_1(J_1)\cap\W_3(J_2)\cap\W_3(J_3)$. In
other words, $(M,J_1,g)$ is an almost Hermitian manifold with
neutral metric of the class $\G_1=\W_1\oplus\W_3\oplus\W_4$
according to the classification in \cite{GrHe}, and $(M,J_\al,g)$,
$\al=2,3$, are almost complex manifolds with Norden metric of the
class $\W_3$ (the co-called quasi K\"ahler manifolds with Norden
metric), according to the classification in \cite{GaBo}. There
these classes are defined as follows:
\begin{gather}
\label{G1}
    \G_1(J_1):\quad
    F_1(x,x,z)=F_1(J_1x,J_1x,z),\\[4pt]
\label{W3}
\W_3(J_\al):\quad \sx F_\al(x,y,z)=0, \quad \al=2,3.
\end{gather}

\begin{thm}\label{thm-G1}
Let $M$ be an almost $(H,G)$-manifold which is a quasi-K\"ahler
manifold with Norden metric regarding $J_2$ and $J_3$. Then it
belongs to the class $\G_1$ with respect to $J_1$.
\end{thm}

\begin{proof}
Suppose  $(M,J_\al,g)$, $\al=2,3$, belongs to the class $\W_3$ of
the quasi-K\"ahler manifolds with Norden metric. Then, according
to \eqref{W3} we have respectively
\begin{gather}
F_2(x,y,z)+F_2(y,x,z)=-F_2(z,x,y),\label{1}\\[4pt]
F_3(x,y,z)+F_3(y,x,z)=-F_3(z,x,y).\label{3}
\end{gather}
The relation $J_2=-J_1\circ J_3$ implies
\begin{equation}
F_2(x,y,z)=F_3(x,y,J_1z)-F_1(x,J_3y,z).\label{2}
\end{equation}
Next, using consecutively \eqref{1}, \eqref{2} and \eqref{3}, we
obtain
\[
\begin{split}
F_2(z,x,y)=&-F_2(x,y,z)-F_2(y,x,z)\\[4pt]
=&-F_3(x,y,J_1z)+F_1(x,J_3y,z)-F_3(y,x,J_1z)+F_1(y,J_3x,z),
\end{split}
\]
\begin{equation}\label{4}
F_2(z,x,y)=F_3(J_1z,x,y)+F_1(x,J_3y,z)+F_1(y,J_3x,z).
\end{equation}
Having in mind $J_3=-J_2\circ J_1$, we express the term
$F_3(J_1z,x,y)$ from the line above as follows
\begin{equation}\label{5}
F_3(J_1z,x,y)=-F_2(J_1z,J_1x,y)-F_1(J_1z,x,J_2y).
\end{equation}
Then, from \eqref{4} and \eqref{5} we have
\begin{equation}\label{6}
\begin{split}
F_2(z,x,y)+F_2(J_1z,J_1x,y)&=F_1(x,J_3y,z)\\[4pt]
&+F_1(y,J_3x,z)-F_1(J_1z,x,J_2y).
\end{split}
\end{equation}
We replace the  substitutions $z\rightarrow J_1z$ and
$x\rightarrow J_1x$ in \eqref{6} and applying  properties
\eqref{34} and \eqref{J123}, we obtain
\begin{equation}\label{7}
\begin{split}
F_2(J_1z,J_1x,y)+F_2(z,x,y)=&-F_1(J_1x,J_2y,z)\\[4pt]
&+F_1(y,J_3x,z)+F_1(z,J_1x,J_2y).
\end{split}
\end{equation}
The difference of \eqref{6} and \eqref{7} is the following
\[
F_1(x,J_3y,z)+F_1(J_1x,J_2y,z)-F_1(J_1z,x,J_2y)-F_1(z,J_1x,J_2y)=0.
\]
Next, we substitute $J_3y$ for $y$ in the last equality and
applying \eqref{34} we obtain
\begin{equation}\label{G1=}
F_1(x,z,y)-F_1(J_1x,J_1z,y)-F_1(J_1z,J_1x,y)+F_1(z,x,y)=0,
\end{equation}
which is equivalent to \eqref{G1}, \ie $(M,J_1,g)$ belongs to the
class $\G_1$.
\end{proof}

\begin{thm}\label{thm-K}
Let $M$ be an almost $(H,G)$-manifold from  the class
$\G_1(J_1)\cap\W_3(J_2)\cap\W_3(J_3)$. If $(M,J_1,g)$ 
belongs to the subclass $\W_0(J_1):\ F_1=0$ of the K\"ahler
manifolds with neutral metric then $(M,H,G)$ is a
pseudo-hyper-K\"ahler manifold.
\end{thm}
\begin{proof}
Suppose $F_1=0$ holds. By virtue of
\eqref{F-prop-1} we have
\begin{equation}\label{F2=F3J1}
    F_2(x,y,z)=F_3(x,J_1 y,z)=F_3(x,y,J_1 z)
\end{equation}
and by substitutions $y \rightarrow J_1y$, $z \rightarrow J_1z$
the following relations are valid
\begin{equation}\label{F2J1=F3J1}
    F_2(x,J_1y,J_1z)=-F_3(x, y,J_1z)=-F_3(x,J_1y,z).
\end{equation}
Comparing \eqref{F2=F3J1} and \eqref{F2J1=F3J1}, we obtain
\begin{equation}\label{F2J1=F2}
    F_2(x,J_1y,J_1z)=-F_2(x,y,z).
\end{equation}
Having in mind the property $F_2(x,J_2y,J_2z)=F_2(x,y,z)$ from
\eqref{34}, equality \eqref{F2J1=F2} implies
\begin{equation}\label{F2J3=F2}
    F_2(x,J_3y,J_3z)=-F_2(x,y,z).
\end{equation}

Now, from \eqref{F2=F3J1} we have
\begin{equation}\label{F2=F3J3}
    F_2(x,y,J_2z)=F_3(x,y,J_3z).
\end{equation}
According to the conditions $\sx F_2(x,y,z)=\sx F_3(x,y,z)=0$ for
a manifold in $\W_3(J_2)\cap\W_3(J_3)$ and since $F_2$ and $F_3$
are symmetric by the second and the third arguments, then by
virtue of \eqref{F2=F3J3} the following equality is valid
\begin{equation}\label{F2J2=F3J3}
    F_2(J_2z,x,y)=F_3(J_3z,x,y).
\end{equation}
Using $F_3(x,J_3y,J_3z)=F_3(x,y,z)$ from \eqref{34}, we obtain
\begin{equation}\label{F2=F2J3}
    F_2(z,x,y)=F_2(z,J_3x,J_3y).
\end{equation}
Equalities \eqref{F2J3=F2} and \eqref{F2=F2J3} imply $F_2=0$ and
therefore $F_3=0$ is also valid, \ie $(M,H,G)\in\K$.
%
%
%
\end{proof}



\subsection{Class $\W_{133}$}

In this subsection we consider a more specialized manifold
$(M,H,G)$ with non-integrable structures $J_\al$ ($\al=1,2,3$),
namely a manifold which is \emph{nearly K\"ahlerian with neutral
metric} regarding $J_1$ and \emph{quasi K\"ahlerian with Norden
metric} regarding $J_2$ and $J_3$. In other words, $(M,H,G)$
belongs to the class $\W_1(J_1)\cap\W_3(J_2)\cap\W_3(J_3)$ (in
short, $\W_{133}$) according to the corresponding classifications
in \cite{Gray} and \cite{GaBo}, where the basic class $\W_1(J_1)$
is defined by
\begin{equation}\label{W1}
    \W_1(J_1):\quad F_1(x,y,z)=-F_1(y,x,z)
\end{equation}
or equivalently $F_1(x,y,z)=\frac{1}{3}\sx
F_1(x,y,z)=\frac{1}{3}\dd g_1(x,y,z)$. Furthermore, the property
$F_1(J_1 x,J_1 y,z)=-F_1(x,y,z)$ holds for a $\W_1(J_1)$-manifold
since $\W_1(J_1)\subset\G_1(J_1)$ and \eqref{G1}.

It is known (\cite{GrHe}) that the class
$\G_1=\W_1\oplus\W_3\oplus\W_4$ of almost Hermitian manifolds
$(M,J_1,g)$ exists in general form when the dimension of $M$ is at
least 6. At dimension 4, $\G_1$ is restricted to its subclass
$\W_4$. Thus, the manifold $(M,H,G)$ belonging to the class
$\W_{133}$ exists when $\dim{M}\geq 8$.

\subsubsection{Properties of the structural (0,3)-tensors in
$\W_{133}$}

According to  \eqref{34} for $\al=1$ and \eqref{W1}, we have
\begin{equation}\label{F1_3-form}
    F_1(x,y,z)=-F_1(y,x,z)=-F_1(x,z,y),
\end{equation}
\ie $F_1$ is a 3-form, and moreover
\begin{gather}
F_1(x,y,z)=- F_1(J_1x,J_1y,z)=- F_1(J_1x,y,J_1z)=- F_1(x,J_1y,J_1z),
\label{F1-hibrid}\\[4pt]
F_1(J_1x,y,z)=F_1(x,J_1y,z)=F_1(x,y,J_1z).\label{F1-J1}
\end{gather}

\begin{lem}\label{lem-Fal}
The structural tensors $F_\al$ ($\al=1,2,3$) of
$(M,H,G)\in\W_{133}$ has the following properties:
\begin{gather}
\begin{split}\label{F1-23}
    F_1(J_2x,J_2y,z)=- F_1(J_3x,J_3y,z),\\[4pt]
    F_1(J_2x,y,J_2z)=- F_1(J_3x,y,J_3z),\\[4pt]
    F_1(x,J_2y,J_2z)=- F_1(x,J_3y,J_3z),
\end{split}
    \\[4pt]
    F_2(x,y,z)=F_2(x,J_1y,J_1z)=F_2(x,J_3y,J_3z),\label{F2-123}\\[4pt]
    F_3(x,y,z)=F_3(x,J_1y,J_1z)=F_3(x,J_2y,J_2z).\label{F3-123}
\end{gather}
\end{lem}
\begin{proof}
The first equality of \eqref{F1-23} follows from the first
equality of \eqref{F1-hibrid} applying $J_2$ to $x$ and $y$. In a
similar way, we obtain the other two properties in \eqref{F1-23}.

We apply $J_3$ to $x$ and $y$ in \eqref{4}.  Then, using
\eqref{34} and \eqref{F1_3-form}, we obtain
$F_2(z,J_3x,J_3y)=F_2(z,x,y)$. The remaining part of
\eqref{F2-123} follows from the latter equality by the
substitutions $x\rightarrow J_2x$, $y \rightarrow J_2y$ and
\eqref{34}.

The proof of \eqref{F3-123} is analogous to the previous one.
\end{proof}

The last lemma imply immediately the following
\begin{lem}\label{lem2-Fal}
The structural tensors $F_\al$ ($\al=1,2,3$) of
$(M,H,G)\in\W_{133}$ has the following properties:
\begin{gather}
\begin{split}\label{F1-23=}
    F_1(J_2x,J_3y,z)=F_1(J_3x,J_2y,z),\\[4pt]
    F_1(J_2x,y,J_3z)=F_1(J_2x,y,J_3z),\\[4pt]
    F_1(x,J_2y,J_3z)=F_1(x,J_2y,J_3z),
\end{split}
    \\[4pt]
    F_2(x,y,J_1z)=-F_2(x,J_1y,z),\qquad F_2(x,y,J_3z)=-F_2(x,J_3y,z),\label{F2-123=}\\[4pt]
    F_3(x,y,J_1z)=-F_3(x,J_1y,z),\qquad F_3(x,y,J_2z)=-F_3(x,J_2y,z).\label{F3-123=}
\end{gather}
\hfill$\Box$
\end{lem}

The latter two lemmas imply the following
\begin{prop}\label{prop-Fal=}
The structural tensors $F_\al$ ($\al=1,2,3$) of a
$\W_{133}$-manifold $(M,H,G)$ have the following properties:
\begin{gather}
    2F_2(x,y,z)=F_1(x,y,J_3 z)-F_1(x,J_3 y,z),\label{2F2}\\[4pt]
    2F_3(x,y,z)=F_1(x,J_2 y,z)-F_1(x,y,J_2 z),\label{2F3}\\[4pt]
    F_2(x,y,z)=-F_3(J_1x,y,z),\label{F2F3J1}\\[4pt]
    F_1(x,y,J_1z)+F_2(x,y,J_2z)+F_3(x,y,J_3z)=0.\label{sumFalJal}
\end{gather}
\end{prop}
\begin{proof}

From \eqref{F-prop-1} we have
\[
\begin{array}{l}
    F_2(x,y,z)=F_3(x,J_1 y,z)+ F_1(x,y,J_3 z),\\[4pt]
    F_2(x,y,z)=-F_1(x,J_3 y,z)+F_3(x,y,J_1 z)
\end{array}
\]
which we sum up, apply \eqref{F3-123=}  so the result is
\eqref{2F2}. Similarly, we obtain \eqref{2F3}.

Next, we replace $z$ with $J_2z$ and $J_3z$ in \eqref{2F2} and
 \eqref{2F3}, respectively. Then, the addition of the resultant
 equalities yields \eqref{sumFalJal},
 because of \eqref{F1-J1}.
\end{proof}

\begin{prop}\label{prop-W133K}
Let $M$ be an almost $(H,G)$-manifold from  the class $\W_{133}$.
If $(M,J_\al,g)$ (for some $\al=1,2,3$) belongs to the subclass
$\W_0(J_\al):\ F_\al=0$ of the corresponding manifolds of K\"ahler
 type then $(M,H,G)$ is a pseudo-hyper-K\"ahler
manifold.
\end{prop}
\begin{proof}
It follows from the interconnections between $F_1$, $F_2$ and
$F_3$ in \propref{prop-Fal=}.
\end{proof}

Having in mind \propref{prop-W133K}, we can restrict our attention
to the class of the the so-called \emph{strict
$\W_{133}$-manifolds} $(M,H,G)$, \ie $(M,J_\al,g)$ does not belong
to $\W_0(J_\al)$ for any $\al=1,2,3$.

\subsubsection{Curvature properties in $\W_{133}$}

Let us recall the following curvature properties on a nearly
K\"ahler manifold $(M,J_1,g)$ known from \cite{Gray}:
\begin{equation}\label{RGray}
    R(x,y,J_1z,J_1w)-R(x,y,z,w)=A_1(x,y,z,w),
\end{equation}
where \(
A_1(x,y,z,w)=g\bigl(\left(\n_xJ_1\right)y,\left(\n_zJ_1\right)w\bigr)
=g\bigl(T(x,y),T(z,w)\bigr)\) and
\begin{equation}\label{RGray2}
    R(J_1x,J_1y,J_1z,J_1w)=R(x,y,z,w).
\end{equation}

\begin{thm}\label{thm-R}
The curvature tensor $R$ on $(M,H,G)\in\W_{133}$ has the following
property with respect to the almost hypercomplex structure $H$:
\begin{equation}\label{R}
\begin{split}
    R(x,y,z,w)&+\sum_{\al=1}^3 R(x,y,J_\al z,J_\al w)\\[4pt]
&=\sum_{\al=1}^3\bigl\{A_\al(x,z,y,w)-A_\al(y,z,x,w)\bigr\},
\end{split}
\end{equation}
where
\(
A_\al(x,y,z,w)=g\bigl(\left(\n_xJ_\al\right)y,\left(\n_zJ_\al\right)w\bigr)
\), $\al=1,2,3$.
\end{thm}
\begin{proof}
We apply the covariant derivative by $\n$ to \eqref{sumFalJal} and
obtain
\[
\sum_{\al=1}^3\bigl\{\left(\n_xF_\al\right)(y,z,J_\al
w)+A_\al(y,z,x,w)\bigr\}=0.
\]
We alternate the last equality with respect to $x$ and $y$.
Further, we apply the following corollary of the Ricci identity on
an $(H,G)$-manifold known from \cite{Man28}:
\begin{equation}\label{Ric-id-J} %
\begin{split}%
\left(\n_x F_\al\right)(y,z,J_\al w)&-\left(\n_y
F_\al\right)(x,z,J_\al
w)\\[4pt]
&=R(x,y,J_\al z,J_\al w)-\ea R(x,y,z,w).
\end{split}
\end{equation}
The result is \eqref{R}.
\end{proof}

\begin{cor}\label{cor-RJ23}
The curvature tensor $R$ on $(M,H,G)\in\W_{133}$ has the
following properties with respect to the almost complex structures $J_2$ and $J_3$, respectively: 
\begin{equation}\label{RJ2}
\begin{split}
    2R(x,y,z,w)+&2R(x,y,J_2 z,J_2 w)\\[4pt]
=&\sum_{\al=1}^3\bigl\{A_\al(x,z,y,w)-A_\al(y,z,x,w)\bigr\}\\[4pt]
&-A_1(x,y,z,w)-A_1(x,y,J_2z,J_2w),
\end{split}
\end{equation}
\begin{equation}\label{RJ3}
\begin{split}
    2R(x,y,z,w)+&2R(x,y,J_3 z,J_3 w)\\[4pt]
=&\sum_{\al=1}^3\bigl\{A_\al(x,z,y,w)-A_\al(y,z,x,w)\bigr\}\\[4pt]
&-A_1(x,y,z,w)-A_1(x,y,J_3z,J_3w).
\end{split}
\end{equation}
\end{cor}
\begin{proof}
From \eqref{RGray} we have
\[
R(x,y,J_1z,J_1w)=R(x,y,z,w)+A_1(x,y,z,w).
\]
After the substitutions $z \rightarrow J_2z$ and $w \rightarrow
J_2w$ in the last equality we obtain
\[
R(x,y,J_3z,J_3w)=R(x,y,J_2z,J_2w)+A_1(x,y,J_2z,J_2w).
\]
We replace the left-hand sides of the last two equalities in
\eqref{sumFalJal} and get \eqref{RJ2}.

In a similar way we obtain \eqref{RJ3}. Moreover, we establish the
property $A_1(x,y,J_2z,J_2w)=-A_1(x,y,J_3z,J_3w)$ because of the
first line of \eqref{F1-23}.
\end{proof}

As an immediate consequence of \coref{cor-RJ23} we obtain for the
scalar curvature and its associated quantities by $J_\al$ the
following relations:
\begin{equation}
    \tau-\tau^{**}_1=\nJ{1},\quad
    \tau+\tau^{**}_\al=-\frac{1}{2}\nJ{\al}, \quad \al=2,3,\label{tau23}
\end{equation}
where $\tau^{**}_\al=g^{ij}g^{ks}R(e_i,e_k,J_\al e_s,J_\al
e_j)$, $\al=1,2,3$.

According to \eqref{F2F3J1} we get directly
\begin{equation}\label{nJ23}
    \nJ{2}=\nJ{3}.
\end{equation}
Using \eqref{2F2} and \eqref{2F3} we obtain
\[
    \nJ{2}=-\frac{1}{2}\nJ{1}
    -\frac{1}{2}g^{ij}g^{ks}g^{pq}F_1(e_i,e_k,e_p)F_1(e_j,J_3e_s,J_3e_q),
\]
\[
    \nJ{3}=-\frac{1}{2}\nJ{1}
    -\frac{1}{2}g^{ij}g^{ks}g^{pq}F_1(e_i,e_k,e_p)F_1(e_j,J_2e_s,J_2e_q).
\]
Since the third property in \eqref{F1-23} holds,  then summing up
the latter two equalities we obtain
\[
\nJ{1}+\nJ{2}+\nJ{3}=0.
\]
Because of \eqref{nJ23} the last equality implies
\begin{equation}\label{nJ123}
    \nJ{1}=-2\nJ{2}=-2\nJ{3}.
\end{equation}
%
From \eqref{tau23} and \eqref{nJ123} we obtain the following
relation
\begin{equation}\label{tau123**}
    3\tau+\tau^{**}_1=-4\tau^{**}_{2}=-4\tau^{**}_{3}.
\end{equation}

\section{Pseudo-HKT-connection on an $(H,G)$-manifold which is
a near\-ly K\"ahler manifold with respect to $J_1$}

\subsection{Constructing the pHKT-connection}

As it is known from \cite{GaMi}, a linear connection $D$ is called
a \emph{natural connection on an almost complex manifold $(M,J)$
with Norden metric} $g$ if the almost complex structure $J$ and
the metric $g$ are parallel with respect to $D$, i.e. $DJ=Dg=0$.

The notion of the \emph{hyper-Hermit\-ian connection} in the
hyper-Hermitian geometry is known. This is a linear connection
such that the three almost complex structures and the Hermitian
metric are parallel regarding this connection.

Following this idea we give

\begin{defn}
A linear connection $D$ is called a \emph{natural connection on an
almost hypercomplex manifold $(M,H)$ with a pseudo-Riemannian
metric $g$} if the almost hypercomplex structure $H=(J_1,J_2,J_3)$
and the metric $g$ are parallel with respect to $D$, i.e.
$DJ_1=DJ_2=DJ_3=Dg=0$.
\end{defn}

As a corollary for a natural connection $D$ on $(M,H,G)$ we have
also $Dg_1=Dg_2=Dg_3=0$, having in mind \eqref{gJ}.

If $T$ is the torsion tensor of $D$, i.e. $T(x,y)=D_x y-D_y x-[x,
y]$, then the corresponding tensor field of type (0,3) is
determined by $T(x,y,z)=g(T(x,y),z)$.


In a similar way to  KT- and HKT-connections  we introduce

\begin{defn}
A natural connection $D$ is called a \emph{pseudo-HKT-connection
on an almost $(H,G)$-manifold} (briefly, \emph{pHKT-connection})
if the torsion tensor $T$ of $D$ is totally skew-symmetric.
\end{defn}

For an almost complex manifold with Hermitian metric $(M,J,g)$, in
\cite{Fri-Iv2} it is proved that there exists a unique
KT-connection if and only if the Nijenhuis tensor $N_{J}(x,y,z):=
g(N_{J}(x,y),z)$ is a 3-form, \ie the manifold belongs to the
class of cocalibrated structures $\G_1=\W_1\oplus\W_3\oplus\W_4$.
KT-connections on nearly K\"ahler manifolds are investigated in
\cite{Nagy07} for instance.


Then, there exists a unique KT-connection $D^1$ for the almost
Hermitian manifold $(M,J_1,g)\in\G_1$ on the considered almost
$(H,G)$-manifold such that
\begin{equation}\label{D1}
    g\left(D^1_xy,z\right)=g\left(\n_xy,z\right)+Q_1(x,y,z),
\end{equation}
where
\begin{equation}\label{Q1}
    Q_1(x,y,z)=\frac{1}{2}\left\{F_1(x,y,J_1z)+F_1(y,z,J_1x)-F_1(z,x,J_1y)\right\}.
\end{equation}

The difference $Q_1$ of $D^1$ and $\n$ is a totally skew-symmetric
tensor because of the relation $T_1=2Q_1$ with the corresponding
torsion $T_1$ of the KT-connection $D^1$. Then, in the case when
$\G_1$ is restricted to $\W_1$ for $J_1$, because of \eqref{W1},
expression \eqref{Q1} assumes the form
\begin{equation}\label{Q1=}
    Q_1(x,y,z)=\frac{1}{2}F_1(x,y,J_1z).
\end{equation}

In \cite{Mek2} it is proved that there  exists a unique
KT-connection on an almost complex manifold with Norden metric
$(M,J,g)$ if and only if
the manifold is quasi-K\"ahlerian with Norden metric or $(M,J,g)$
belongs to the class $\W_3$.

Then, there exists a unique KT-connection $D^\al$ ($\al=2,3$,
respectively) for the almost complex manifold with Norden metric
$(M,J_\al,g)\in\W_3(J_\al)$
on the considered almost
$(H,G)$-manifold such that
\begin{equation}\label{D23}
    g\left(D^\al_xy,z\right)=g\left(\n_xy,z\right)+Q_\al(x,y,z),
\end{equation}
where
\begin{equation}\label{Q23}
    Q_\al(x,y,z)=-\frac{1}{4}\sx F_\al(x,y,J_\al z).
\end{equation}

Let us construct a linear connection $D$, using the KT-connections
$D^1$, $D^2$ and $D^3$, on an almost $(H,G)$-manifold from the
class $\G_1(J_1)\cap\W_3(J_2)\cap\W_3(J_3)$ as follows
\begin{equation}\label{D}
    g\left(D_xy,z\right)=g\left(\n_xy,z\right)+Q(x,y,z),
\end{equation}
where
\begin{equation}\label{Q}
    Q(x,y,z)=-\frac{1}{2}Q_1(x,y,z)+Q_2(x,y,z)+Q_3(x,y,z).
\end{equation}




We apply \eqref{sumFalJal}, \eqref{34}, \eqref{W1} and
\eqref{F1-J1} in \eqref{Q} and obtain
\[
Q(x,y,z)=\frac{1}{2}F_1(x,y,J_1z).
\]

Then, the linear connection $D$, introduced by \eqref{D} and
\eqref{Q}, has the following definitional equality
\begin{equation}\label{DFW1}
    g\left(D_xy,z\right)=g\left(\n_xy,z\right)+\frac{1}{2}F_1\left(x,y,J_1z\right)
\end{equation}
or equivalently
\begin{equation}\label{DW1}
    D_xy=\n_xy-\frac{1}{2}J_1\left(\n_xJ_1\right)y.
\end{equation}

We verify directly that $D$ is a natural connection on $(M,H,G)$
and therefore, because of the properties of $F_1$ and the
uniqueness of $D^1$, $D^2$ and $D^3$, we obtain
\begin{prop}\label{prop-HKT}
The linear connection $D$ defined by \eqref{DW1} is a unique
pHKT-connection on an almost $(H,G)$-manifold from the class
$\W_{133}$.\hfill$\Box$
\end{prop}

Let us remark that the pHKT-connection $D$ on an almost
$(H,G)$-manifold coincides with the known KT-connection $D^1$ on
nearly K\"ahler manifolds studied in \cite{Fri-Iv2}, \cite{Nagy07}
for the Riemannian case and \cite{CorSch} for the
pseudo-Riemannian case.

\subsection{Characteristics of the constructed pHKT-connection}
\subsubsection{$D$-parallel torsion of $D$}

Since $D$ coincides with $D^1$, then the pHKT-connection $D$
defined by \eqref{DW1} on an almost $(H,G)$-manifold from
$\W_{133}$ is $D$-parallel, $\dd T=0$, and henceforth $T$ is
coclosed, $\delta T=0$ (\cite{Fri-Iv2}, \cite{Kir},
\cite{BelMor}).

\subsubsection{Curvature tensor of $D$} Let us consider the
curvature tensor $K$ of the connection $D$, \ie $K(x,y)z=[D_x,
D_y] z - D_{\left[x,y\right]} z$ and
$K(x,y,z,w)=g\left(K(x,y)z,w\right)$. The relation between the
connections $D$ and $\n$ generates the corresponding relation
between their curvature tensors $K$ and $R$.

\begin{prop}\label{prop-K}
Let $(M,H,G)$ be an almost $(H,G)$-manifold from $\W_{133}$ and
$D$ be the pHKT-connection defined by \eqref{DW1}. Then the
cur\-va\-ture tensors $K$ of $D$ and $R$ of $\n$ has the following
relation
\begin{equation}\label{KR}
    K(x,y,z,w)=R(x,y,z,w)+\frac{1}{4} A_1(x,y,z,w)+\frac{1}{4}\sx A_1(x,y,z,w).
\end{equation}
\hfill$\Box$
\end{prop}

Let us remark that such a relation is given in \cite{Gray76},
\cite{Nagy07} and \cite{Fri-Iv2}.

\propref{prop-K} implies directly that $\rho^D$ is symmetric and
$J_1$-invariant, which is also known from \cite{Fri-Iv2} and
\cite{Nagy07} for the Riemannian case. Moreover, for the scalar
curvatures we have
\begin{equation}\label{tauD}
\tau^D=\tau^\n-\frac{1}{4}\norm{T},
\end{equation}
where $\norm{T}=\nJ{1}$.

Further, we give some necessary and sufficient conditions for
flatness of $D$ in the following
\begin{prop}\label{prop-sA}
Let $(M,H,G)$ be an almost $(H,G)$-manifold from $\W_{133}$ and
$D$ be the parallel pHKT-connection defined by \eqref{DW1}. Then
the following characteristics of this connection are equivalent:
\begin{enumerate}\renewcommand{\labelenumi}{(\roman{enumi})}
    \item $D$ is strong;
    \item $D$ has a $\n$-parallel torsion;
    \item $D$ is flat.
\end{enumerate}
\end{prop}

\begin{proof}
Using $DT=0$ we obtain
\begin{equation}\label{nTA}
    \left(\n_xT\right)(y,z,w)=\frac{1}{2}\sx A_1(x,y,z,w).
\end{equation}
Then we compute the exterior derivative of the 3-form $T$ as
\begin{equation}\label{dT}
\dd T(x,y,z,w)=2\sx A_1(x,y,z,w).
\end{equation}
Therefore, the pHKT-connection is strong, \ie $\dd T=0$, if and
only if the identity  $\sx A_1(x,y,z,w)=0$ is valid, which is
equivalent to $\n T=0$, according to \eqref{nTA}. Thus, (i) is
equivalent to (ii).

The curvature tensor $K$ satisfies the properties of the first
line of \eqref{curv} and the K\"ahler-like property \eqref{L-kel},
since $DJ_\al=0$, $\al=1,2,3$. Then, $K$ is a K\"ahler-like
 tensor if and only if the identity $\sx K(x,y,z,w)=0$ holds.

On the other hand, the first Bianchi identity for $K$ with torsion
$T$
\[
\sx K(x,y,z,w)=\sx\bigl\{
T\bigl(T(x,y),z,w\bigr)+\left(D_xT\right)(y,z,w)\bigr\}
\]
takes the following form, since $T$ is $D$-parallel:
\begin{equation}\label{KA}
\sx K(x,y,z,w)=\sx A_1(x,y,z,w).
\end{equation}

Therefore, according to \eqref{dT} and \eqref{KA}, we have that
$D$ is a strong pHKT-connection if and only if the curvature
tensor $K$ of $D$ is a K\"ahler-like tensor. Then, because of
\propref{th-0}, $K$ is zero, \ie the equivalence
(i)$\Leftrightarrow$(iii) is valid. This completes the proof.
\end{proof}


\begin{prop}\label{prop-flat}
Let $(M,H,G)$ be an almost $(H,G)$-manifold from $\W_{133}$ and
$D$ be the parallel pHKT-connection defined by \eqref{DW1}. If $D$
is flat or strong then $(M,H,G)$ is $\n$-flat,
isotropic-hyper-K\"ahlerian and the torsion of $D$ is isotropic.
\end{prop}

\begin{proof}
According to \propref{prop-sA}, if $(M,H,G)$ is $D$-flat then we
obtain the following form of the curvature tensor $R$ of the
Levi-Civita connection $\n$, applying \eqref{KR}:
\begin{equation}\label{RA}
    R(x,y,z,w)=
    -\frac{1}{4}A_1(x,y,z,w).
\end{equation}
Substituting $J_1z$ and $J_1w$ for $z$ and $w$ in \eqref{RA},
respectively, and using the properties
$A_1(x,y,J_1z,J_1w)=-A_1(x,y,z,w)$ and \eqref{RGray}
 for a nearly K\"ahler
manifold $(M,J_1,g)$, we obtain $A_1=0$. Therefore, according to
\eqref{RA}, $R=0$ and then the connections $\n$ and $D$ are both
flat. Moreover, having in mind  \eqref{nJ123} and \eqref{tauD},
$(M,H,G)$ is isotropic-hyper-K\"ahlerian and  $T$ is isotropic.
\end{proof}

\begin{cor}
Let $(M,H,G)$ be an almost $(H,G)$-manifold from $\W_{133}$. If it
is not $D$-flat then it does not admit a strong pHKT-connection,
\ie any non-$D$-flat $(M,H,G)$ is a weak
pHKT-manifold.\hfill$\Box$
\end{cor}

It is known from \cite{CorSch} that there are no strict flat
nearly pseudo-K\"ahler manifolds of dimension less than 12.


\end{document}